\documentclass[12pt]{article}

\usepackage[left=2cm,right=2cm,top=4cm,bottom=4cm]{geometry}
\usepackage[T1]{fontenc}
\usepackage[latin1]{inputenc}
\usepackage{lmodern}
\usepackage[nottoc, notlof, notlot]{tocbibind}
\usepackage{amsmath}
\usepackage{amsthm}
\usepackage{amssymb}
\usepackage{ stmaryrd }
\usepackage{graphicx}
\usepackage{tikz} 
\usepackage{dsfont}
\usepackage{xspace}
\usepackage{enumitem}
\usepackage{ upgreek }
\usetikzlibrary{matrix,arrows} 
\usepackage[colorlinks,urlcolor=blue,linkcolor=blue,citecolor=brown]{hyperref}

\theoremstyle{plain}
\newtheorem{thm}{Theorem}

\theoremstyle{remark}

\theoremstyle{plain}

\theoremstyle{definition}

\theoremstyle{plain}

\theoremstyle{plain}

\theoremstyle{plain}

\theoremstyle{remark}

\begin{document}

\title{Note on the spectrum of classical and uniform exponents of Diophantine approximation}  \author{Antoine MARNAT\footnote{supported by the Austrian Science Fund (FWF), Project F5510-N26, and FWF START project Y-901 and EPSRC Programme Grant EP/J018260/1} \\
 \href{mailto:antoine.marnat@york.ac.uk}{antoine.marnat@york.ac.uk}}
 
 \date{}
\maketitle

\abstract{  Using the Parametric Geometry of Numbers introduced recently by W.M. Schmidt and L. Summerer \cite{SchSu,SchSu2} and results by D. Roy \cite{RoyParam,RoySpec}, we establish that  the $2n$ exponents of Diophantine approximation in dimension $n\geq3$ are algebraically independent.}

 \section{Introduction}
 
Throughout this paper, the integer $n\geq1$ denotes the dimension of the ambient space $\mathbb{R}^n$ endowed with its Euclidean norm and $\boldsymbol{\theta}=(\theta_1, \ldots , \theta_n)$ denotes an $n$-tuple of real numbers such that $1,\theta_1, \ldots , \theta_n$ are $\mathbb{Q}$-linearly independent.\\
 
Let $d$ be an integer with $0\leq d \leq n-1$. We define the exponent ${\omega}_{d}(\boldsymbol{\theta})$  (resp. the uniform exponent $\hat{\omega}_{d}(\boldsymbol{\theta})$) as the supremum of the real numbers $\omega$ for which there exist  rational affine subspaces $L \subset \mathbb{R}^{n}$ such that 
 \[ \dim(L)=d \; ,  \; H(L)\leq H \; \textrm{ and } \; H(L)d(\boldsymbol{\theta},L) \leq H^{-\omega}  \]
  for arbitrarily large real numbers $H$  (resp. for every sufficiently large real number $H$). Here $H(L)$ denotes the exponential height of $L$ (see \cite{SchH} for more details), and $d(\boldsymbol{\theta},L)=\min_{P\in L} d(\boldsymbol{\theta},P)$ is the minimal distance between $\boldsymbol{\theta}$ and a point of $L$. Note that this definition is independent of the choice of a norm on $\mathbb{R}^n$.\\

These exponents were introduced originally by M. Laurent \cite{MLwd}. They interpolate between the classical exponents $\omega(\boldsymbol{\theta})=\omega_{n-1}(\boldsymbol{\theta})$ and $\lambda(\boldsymbol{\theta})=\omega_0(\boldsymbol{\theta})$ (resp. $\hat{\omega}(\boldsymbol{\theta})=\hat{\omega}_{n-1}(\boldsymbol{\theta})$ and $\hat{\lambda}(\boldsymbol{\theta})=\hat{\omega}_0(\boldsymbol{\theta})$) that were  introduced by A. Khintchine \cite{Khin2,Khin1}, V. Jarn\'ik \cite{JAR} and Y. Bugeaud and M. Laurent \cite{BugLau2,BugLau}.\\
 
 We have the relations
 \[ \omega_{0}(\boldsymbol{\theta}) \leq \omega_{1}(\boldsymbol{\theta}) \leq \cdots \leq \omega_{n-1}(\boldsymbol{\theta}) ,  \]
 \[\hat{\omega}_{0}(\boldsymbol{\theta})\leq \hat{\omega}_{1}(\boldsymbol{\theta}) \leq \cdots \leq \hat{\omega}_{n-1}(\boldsymbol{\theta}),\]
 and Minkowski's First Convex Body Theorem \cite{Mink} and Mahler's compound convex bodies theory provide the lower bounds
 
 \begin{equation*}
 \omega_{d}(\boldsymbol{\theta}) \geq \hat{\omega}_{d}(\boldsymbol{\theta}) \geq \frac{d+1}{n-d}, \; \textrm{  for  } 0 \leq d \leq n-1 .
 \end{equation*}
 
 These $2n$ exponents happen to be related as was first noticed by Khinchin with his transference theorem \cite{Khin1}. We use the following notion of spectrum to study more general transfers. Given $k$ exponents $e_{1}, \ldots , e_{k}$, we define the \emph{spectrum} of the exponents $(e_{1}, \ldots , e_{k})$ as the subset of $\mathbb{R}^{k}$ described by all $k$-uples $ (e_{1}(\boldsymbol{\theta}) , \ldots , e_{k}(\boldsymbol{\theta}) )$ as $\boldsymbol{\theta}$ runs through all points $\boldsymbol{\theta}=(\theta_1, \ldots , \theta_n) \in \mathbb{R}^n$ such that $1,\theta_1, \ldots , \theta_n$ are $\mathbb{Q}$-linearly independent.\\
   
   In \cite{Moi2}, the author proved the following theorem.
   
\begin{thm}\label{ThAIF}
For every integer $n\geq3$, the $n$ uniform exponents $\hat{\omega}_0, \ldots, \hat{\omega}_{n-1}$ are algebraically independent.
\end{thm}

   Using the same construction, it is even possible to show that for every integer $n\geq3$, the spectrum of $\hat{\omega}_0, \ldots, \hat{\omega}_{n-1}$ is a subset of $\mathbb{R}^{n}$ with non empty interior. In this paper, we extend this result as follows. 
   
   \begin{thm}\label{IndepTh}
  For every integer $n\geq3$, the $2n$ exponents $\hat{\omega}_0, \ldots, \hat{\omega}_{n-1},\omega_{0}, \ldots , \omega_{n-1}$ are algebraically independent.
   \end{thm}
   
   In dimension $n=2$, the spectrum is fully described in \cite{ML}:
   \begin{thm}[Laurent, 2009]
   In dimension $2$, the spectrum of the four exponents $\omega_0, \omega_1, \hat{\omega}_0, \hat{\omega}_1$ is described by the inequalities
   \begin{eqnarray}\notag
     \hat{\omega}_0 + 1/\hat{\omega}_1 =1, &
   2 \leq \hat{\omega}_1 \leq +\infty, &   \cfrac{\omega_1(\hat{\omega}_1-1)}{\omega_1 + \hat{\omega}_1} \leq \omega_0 \leq \cfrac{\omega_1 - \hat{\omega}_1 +1}{\hat{\omega}_1}.
   \end{eqnarray}
   When $\hat{\omega}_1 < \omega_1 = +\infty$ we have to understand these relations as $\hat{\omega}_1 -1 \leq \omega_0 \leq +\infty $ and when $\hat{\omega}_1= +\infty$, the set of constraints should be interpreted as ${\omega}_0=\omega_1=+\infty$ and $\hat{\omega}_0=1$.
   \end{thm}

The first equality, relating the two uniform exponents, is known as Jarn\'ik's relation \cite{JAR} and breaks the algebraic independence. Note that this sharpens previously mentioned relations. In dimension $n=1$ the uniform exponent is always equal to $1$.\\

We refer the reader to \cite[§2]{Moi2} for the notation and the presentation of the parametric geometry of numbers, main tool of the proof. We mainly use the notation introduced by D. Roy in \cite{RoyParam,RoySpec} which is essentially dual to the one of W. M. Schmidt and L. Summerer \cite{SchSu,SchSu2}.\\


\section{Proof of the main Theorem \ref{IndepTh} }\label{dem}

To prove Theorem \ref{IndepTh}, we place ourselves in the context of parametric geometry of numbers. We fully use Roy's theorem \cite[Theorem 5]{Moi2} that reduces the study of spectra of Diophantine approximation to the study of the combinatorial properties of generalized $n$-systems. We construct explicitly a family of generalized $(n+1)$-systems with $2n$ parameters, which provides the algebraic independence in the spectrum via Roy's theorem.\\

We fix the dimension $n\geq3$. Consider any family of positive parameters 
\[A_{1}=A_{2}< A_{3} < \cdots < A_{n+1} \; , \; B_{2}< B_{3} < \cdots < B_{n} \; , \; C \; , \; D\]
satisfying the following properties for $2\leq k \leq n$:
\begin{equation}
\begin{split}
 A_{1} + A_{2}+ \cdots + A_{n+1}&=1\; , \; B_{2} < D < CA_{2}, \\
 A_{k+1} < B_{k} < A_{k+2} \; &, \; B_{k}< CA_{k},
\end{split}\label{cond1}
\end{equation}
where $A_{n+2}=\infty$.\\

\begin{figure}[!h] 
 \begin{center}
 \begin{tikzpicture}[scale=0.4]
 
\fill (5,4) circle[radius=4pt];
\fill (13,10) circle[radius=4pt];
\fill (16,12) circle[radius=4pt];
\fill (20,16) circle[radius=4pt];
\fill (1,1) circle[radius=4pt];
\fill (24,17) circle[radius=4pt];
\fill (28,10)circle[radius=4pt];
\fill (31,1) circle[radius=4pt]; 

\draw (9,8) node {$\cdots$};
\draw (9,15) node {$\cdots$};
\draw (9,3) node {$\cdots$};

\draw (18,8) node {$\cdots$};
\draw (18,15) node {$\cdots$};
\draw (18,3) node {$\cdots$};

\draw (29.5,8) node {$\cdots$};
\draw (29.5,15) node {$\cdots$};
\draw (29.5,3) node {$\cdots$};

\draw (26.5,8) node {$\cdots$};
\draw (26.5,15) node {$\cdots$};
\draw (26.5,3) node {$\cdots$};
 
\draw[black, semithick] (3,3)--(0,0) node [above,black] {};
\draw[black, semithick] (8,1)--(0,1) node [left,black] {$A_{2}$};
\draw[black, semithick] (3,3)--(0,3) node [left,black] {$A_{3}$};
\draw[black, semithick] (3,3)--(5,4) node [left,black] {};
\draw[dashed,black] (5,4)--(0,4) node [left,black] {$B_{2}$};
\draw[dashed,black] (16,12)--(0,12) node [left,black] {$B_{k}$};
\draw[dashed, black] (13,10)--(0,10)node [left,black] {$B_{k-1}$};
\draw[black, semithick] (5,4)--(8,4) node [left,black] {};
\draw[black, semithick] (5,4)--(6,5) node [left,black] {};
\draw[black, semithick] (6,5)--(0,5) node [left,black] {$A_{4}$};
\draw[black, semithick] (6,5)--(0,5) node [midway,above,black] {$P_{4}$};
\draw[black, semithick] (8,9)--(0,9) node [left,black] {$A_{k}$};
\draw[black, semithick] (8,9)--(0,9) node [above,midway,black] {$P_{k}$};
\draw[black, semithick] (8,11)--(0,11) node [left,black] {$A_{k+1}$};
\draw[black, semithick] (8,11)--(0,11) node [midway,above,black] {$P_{k+1}$};
\draw[black, semithick] (8,16)--(0,16) node [left,black] {$A_{n+1}$};
\draw[black, semithick] (8,16)--(0,16) node [above,midway,black] {$P_{n+1}$};
\draw[black, semithick] (6,5)--(8,6) node [left,black] {};

\draw[black, semithick] (10,8)--(11,9) node [left,black] {};
\draw[black, semithick] (11,9)--(13,10) node[midway,above,sloped,black] { } ;
\draw[black, semithick] (11,9)--(10,9) node [left,black] {};
\draw[black, semithick] (13,10)--(17,10) node [left,black] {};
\draw[black, semithick] (13,10)--(14,11) node[midway,above,sloped,black] {} ;

\draw[black, semithick] (14,11)--(10,11) node [left,black] {};
\draw[black, semithick] (14,11)--(16,12)  node[midway,above,sloped,black] { } ;
\draw[black, semithick] (16,12)--(17,12) node [left,black] {};
\draw[black, semithick] (16,12)--(17,13) node [left,black] {};
\draw[black, semithick] (10,1)--(17,1);
\draw[black, semithick] (10,4)--(17,4);
\draw[black, semithick] (10,16)--(17,16);

\draw[black, semithick] (19,15)--(20,16) node [left,black] {};
\draw[black, semithick] (20,16)--(19,16) node [left,black] {};
\draw[black, semithick] (20,16)--(22,17) node [left,black] {};
\draw[black, semithick] (22,17)--(24,17) node [right,black] {};
\draw[black, semithick] (24,17)--(25,18) ;
\draw[black, semithick] (22,17)--(24,19) ;
\draw[black, semithick] (25,18)--(26,18) ; %
\draw[black, semithick] (27,18)--(29,18) ; %
\draw[black, semithick] (30,18)--(38,18) node [right,black] {$CA_{n}$};
\draw[black, semithick] (30,18)--(38,18) node [below,midway,black] {$P_{n}$};

\draw[black, semithick] (24,19)--(26,19) node [right,black] {};
\draw[black, semithick] (27,19)--(29,19) node [right,black] {};
\draw[black, semithick] (30,19)--(38,19) node [right,black] {$CA_{n+1}$};
\draw[black, semithick] (30,19)--(38,19) node [above,midway,black] {$P_{n+1}$};

\draw[black, semithick] (19,1)--(26,1) node [midway,above,black] {$P_{1}$};
\draw[black, semithick] (27,1)--(29,1);
\draw[black, semithick] (30,1)--(31,1);
\draw[black, semithick] (31,1)--(36,5);
\draw[dashed,black] (36,5)--(38,5) node [right,black] {$D$};
\draw[black, semithick] (36,5)--(38,6)node [right,black] {$CA_{2}$};
\draw[black, semithick] (36,5)--(31,5)node [above,midway,black] {$P_{2}$};

\draw[black, semithick] (19,4)--(26,4) node [above,midway,black] {$P_{2}$};
\draw[black, semithick] (27,4)--(29,4) node [right,black] {};
\draw[black, semithick] (30,4)--(31,5) node [right,black] {};
\draw[black, semithick] (31,5)--(36,5) node [right,black] {};

\draw[black, semithick] (19,10)--(26,10) node [above,midway,black] {$P_{k-1}$};
\draw[black, semithick] (27,10)--(28,10) node [right,black] {};
\draw[black, semithick] (28,10)--(29,11) node [right,black] {};
\draw[black, semithick] (30,11)--(38,11) node [right,black] {$CA_{k-1}$};
\draw[black, semithick] (30,11)--(38,11) node [above,midway,black] {$P_{k-1}$};

\draw[black, semithick] (19,12)--(26,12) node [above,midway,black] {$P_{k}$};
\draw[black, semithick] (27,12)--(28,13);
\draw[black, semithick] (28,13)--(29,13);
\draw[black, semithick] (30,13)--(38,13)node [right,black] {$CA_{k}$};
\draw[black, semithick] (30,13)--(38,13)node [above,midway,black] {$P_{k}$};

\draw[dashed, black] (1,18)--(1,0) node [below] { $1$};
\draw[dashed, black] (38,19)--(38,0) node [below] { $C$};
\draw[dashed, black] (22,17)--(22,-1) node [below] { $\delta_{n,2}=\mu_{n+1}$};
\draw[dashed, black] (11,9)--(11,0) node [below] { $\delta_{k-1,1}$};
\draw[dashed, black] (13,10)--(13,-1) node [below] {$\delta_{k-1,2}$};
\draw[dashed, black] (14,11)--(14,0) node [below] {$\delta_{k,1}$ };
\draw[dashed, black] (16,12)--(16,-1) node [below] { $\delta_{k,2}$};
\draw[dashed, black] (3,3)--(3,0) node [below] { $\delta_{2,1}$};
\draw[dashed, black] (24,19)--(24,0) node [below] { $\mu_{n}$};
\draw[dashed, black] (31,19)--(31,0) node [below] { $\mu_{1}$};
\draw[dashed, black] (28,19)--(28,0) node [below] { $\mu_{k-1}$};
\draw[dashed,black] (36,5)--(36,0) node [below] {$\mu_0$};

 \end{tikzpicture}
 \end{center}
 \caption{Pattern of the combined graph of $\boldsymbol{P}$ on the fundamental interval $[1,C]$ }\label{fig}
 \end{figure}
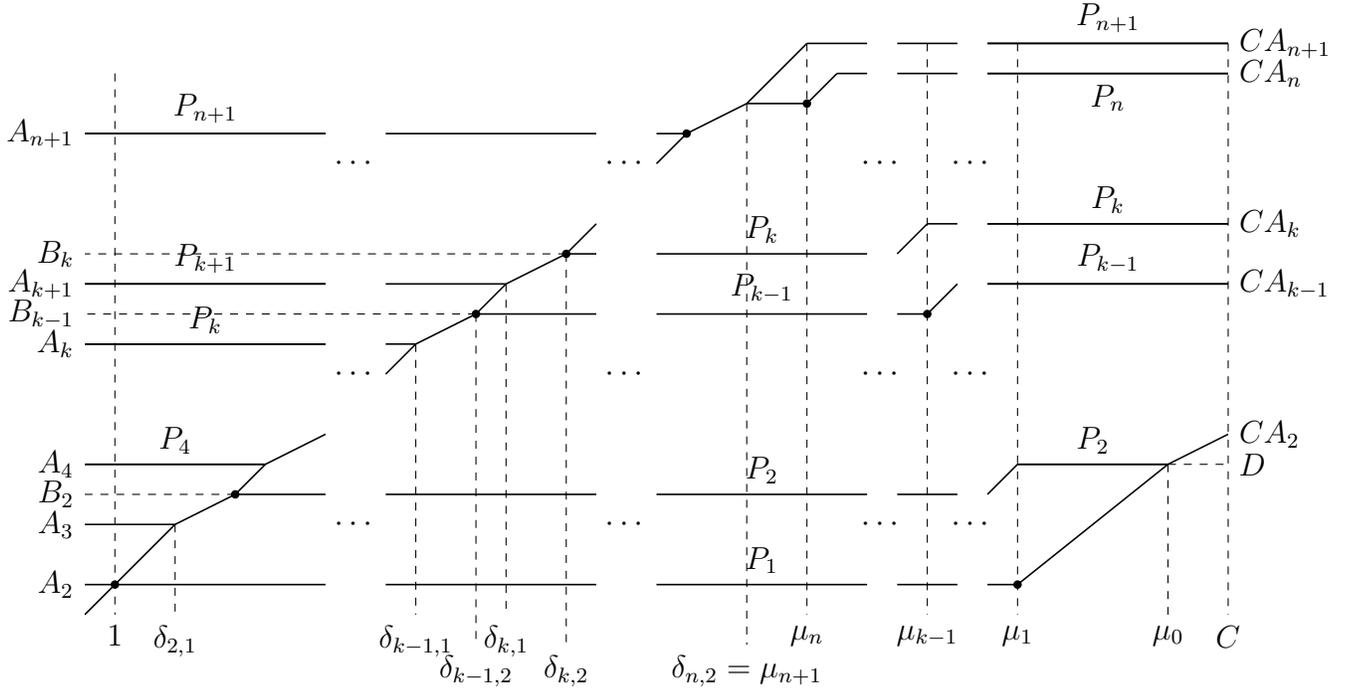

We consider the generalized $(n+1)$-system $\boldsymbol{P}$ on the interval $[1,C]$ depending on the previous parameters whose combined graph is given below by Figure \ref{fig}, where

\[ P_k(1) = A_k \textrm{ and } P_{k}(C)=CA_{k} \textrm{ for } 1\leq k \leq n+1.\]

Conditions \eqref{cond1} are consistent with the graph. On each interval between two consecutive division points, there is only one line segment with non zero slope. This line segment has slope $1$ on the intervals $[1,\delta_{2,1}]$, $[\delta_{k-1,2},\delta_{k,1}]$ for $3\leq k \leq n$, and $[\mu_k,\mu_{k-1}]$ for $n \geq k \geq 1$, and has slope $1/2$ on the interval $[\mu_0,C]$ and $[\delta_{k,1},\delta_{k,2}]$ for $3\leq k \leq n$ , where the two components $P_{k}$ and $P_{k+1}$ coincide. We have $ 3n+1$ division points $1$, $C$, $\delta_{k,1}$ and $\delta_{k,2}$ for $2\leq k \leq n$ and $\mu_l$ for $n+1 \geq l \geq 0$. They are all ordinary division points except $\mu_k$ for $1 \leq k \leq n$ which are switch points.\\

The points which will be most relevant for the proof are labeled with black dots. Note that from $1$ to $\delta_{n,2}$, the combined graph is the same as in \cite[\S5]{Moi2}.\\

We extend $\boldsymbol{P}$ to the interval $[1,\infty)$ by self-similarity. This means, $\boldsymbol{P}(q)=C^{m}\boldsymbol{P}(C^{-m}q)$ for all integers $m$. In view of the value of $\boldsymbol{P}$ and its derivative at $1$ and $C$, one sees that the extension provides a generalized $(n+1)$-system on $[1,\infty)$.\\

The relation between exponents and $n$-systems \cite[Proposition 1]{Moi2} suggests to define $2n$ quantities $W_{n-1}, \ldots , W_{0},\hat{W}_{n-1}, \ldots , \hat{W}_{0}$ by
 \begin{align*}
\cfrac{1}{1+\hat{W}_{n-k}} &:=  \limsup_{q\to+\infty} \cfrac{P_{1}(q) + \cdots + P_k(q)}{q} \textrm{ for } 1\leq k \leq n,\\
 \cfrac{1}{1+W_{n-k}} &:= \liminf_{q\to+\infty} \cfrac{P_{1}(q) + \cdots + P_k(q)}{q}\textrm{ for } 1\leq k \leq n.
   \end{align*}
   Indeed with this setting, Roy's Theorem provides the existence of a point $\boldsymbol{\theta}$ in $\mathbb{R}^{n}$ such that $\hat{\omega}_{k}(\boldsymbol{\theta})= \hat{W}_{k}$ and $\omega_{k}(\boldsymbol{\theta}) = W_{k}$ for every $0\leq k \leq n-1$.\\
  
  Here, self-similarity ensures that the $\limsup$ (resp. $\liminf$) is in fact the maximum (resp. the minimum) on the interval $[1,C[$. Note that for $1\leq k \leq n$, the function $P_1 + \cdots + P_k$ has slope $1$ on the intervals $[1,\delta_{k,1}]$ and $[\mu_{k},C[$, slope $1/2$ on the interval $[\delta_{k,1},\delta_{k,2}]$ and is constant on the interval $[\delta_{k,2},\mu_{k}]$. Therefore the minimum of the function $q \mapsto q^{-1}(P_{1}(q) + \cdots + P_{k}(q))$ is reached at $\mu_{k}$ and its maximum is reached either at $\delta_{k,1}$ or at $\delta_{k,2}$, when slope changes from $1$ to $1/2$ or from $1/2$ to $0$. Namely, the maximum is reached at $\delta_{k,1}$ if 
\begin{equation}\label{delta}
\cfrac{P_{1}(\delta_{k,1}) + \cdots +P_k(\delta_{k,1}) }{\delta_{k,1}}\geq\cfrac{1}{2}
\end{equation}
and at $\delta_{k,2}$ if the lefthand side is $\leq 1/2$. We deduce that for $1\leq k \leq n$,
\begin{eqnarray*}
\hat{W}_{n-k}& =& \cfrac{P_{k+1}(q) + \cdots + P_{n+1}(q)}{P_{1}(q) + \cdots + P_{k}(q)} \; \textrm{ where } q=\left\{ \begin{array}{cl} \delta_{k,1} &\textrm{ if \eqref{delta} is satisfied }\\ \delta_{k,2} & \textrm{ otherwise} \end{array}\right. ,\\
W_{n-k} &=& \cfrac{P_{k+1}(\mu_{k}) + \cdots + P_{n+1}(\mu_{k})}{P_{1}(\mu_{k}) + \cdots + P_{k}(\mu_{k})}.
\end{eqnarray*}

   It is easy to check that the parameters
  \begin{equation}
  \begin{split}
C=3, \, A_{1}=A_{2}=2^{-n}\, , A_{k}&=2^{-n+k-2}\; \textrm{ for } 3 \leq k \leq n+1\\
D=\cfrac{11}{8}2^{-n+1} \, , B_{k}&=\cfrac{5}{4}2^{-n+k-1} \; \textrm{ for } 2 \leq k \leq n
\end{split}\label{p}
 \end{equation}
satisfy the conditions \eqref{cond1}.  For this choice of parameters, the lefthand side of inequality \eqref{delta} is $>1/2$ for $1\leq k \leq n-1$ and  $<1/2$ for $k=n$. This property remains true for $(C,A_{2}, \ldots , A_{n},D,B_{2}, B_{3}, \ldots , B_{n})$ in an open neighborhood of the point \[(3,2^{-n}, \ldots , 2^{-2}, \cfrac{11}{8}2^{-n+1},\cfrac{5}{2}2^{-n}, \ldots , \cfrac{5}{2}2^{-2} )\] provided that we set $A_{1}=A_{2}$ and $A_{n+1}=1-(A_{1}+ \cdots + A_{n})$. In this neighborhood, the quantities $W_{0}, \ldots , W_{n-1}, \hat{W}_{0}, \ldots , \hat{W}_{n-1}$ are given by the following rational fractions in $\mathbb{Q}(C,A_{2}, \ldots , A_{n},D,B_{2}, B_{3}, \ldots , B_{n})$ :

\[\begin{array}{ll}
\hat{W}_{n-1}= \cfrac{1}{A_{2}} -1  , & \hat{W}_{0}  = \cfrac{1- (2A_{2} + A_{3} + A_{4} + \cdots +A_{n}) }{A_2 + (B_2 + \cdots + B_{n-1}) },\\
\hat{W}_{n-k} =  \cfrac{1- (2A_2 + A_{3} +A_{4}+ \cdots + A_{k+1})+B_{k}}{A_2 + (B_2 + \cdots + B_{k})}   & \textrm{ for } 2 \leq k \leq n-1,\\
 W_{n-k} = \cfrac{C(1-(2A_{2} + A_{3} + A_{4} + \cdots + A_{k}))}{A_{2}+B_{2}+ \cdots + B_{k}}  & \textrm{ for } 2 \leq k \leq n,\\
 W_{n-1} = \cfrac{D + C(1-2A_{2}) }{A_{2}}.
 \end{array}\]

Since $W_{0}, \ldots , W_{n-1}, \hat{W}_{0}, \ldots , \hat{W}_{n-1}$ come from a generalized $(n+1)$-system $\boldsymbol{P}$, Roy's Theorem provides the existence of a point $\boldsymbol{\theta}$ in $\mathbb{R}^{n}$ such that $\hat{\omega}_{k}(\boldsymbol{\theta})= \hat{W}_{k}$ and $\omega_{k}(\boldsymbol{\theta}) = W_{k}$ for every $0\leq k \leq n-1$. Therefore, to prove Theorem \ref{IndepTh} it is sufficient to show that the rational fractions $W_{0}, \ldots , W_{n-1},\hat{W}_{0}, \ldots , \hat{W}_{n-1}\in \mathbb{Q}(C,A_{2}, A_{3}, \ldots , A_{n},D,B_{2}, B_{3}, \ldots , B_{n})$ are algebraically independent.\\

First, note that only $W_{n-1}$ depends on $D$ and $\hat{W}_{n-1}$ only depends on $A_{2}$. Therefore, it is enough to prove that the $2n-2$ other rational fractions are algebraically independent over $\mathbb{Q}(A_2)$. For the calculation, it is convenient to successively make the following two changes of variables. First, we set
\begin{eqnarray*}
 M_{k} &:=& 1- \sum_{i=1}^{k} A_{i} \; \textrm{ for } 2 \leq k \leq n+1, \\
 N_{k} &:=& A_{1} + \sum_{i=2}^{k}B_{i} \; \textrm{ for } 1\leq k \leq n .
 \end{eqnarray*}
Note that $M_{n+1}=0$ and $N_1=A_1$. We get the formulae
\begin{eqnarray*}
\hat{W}_{0} &= &\cfrac{M_{n}}{N_{n-1}},\\
W_{n-k}&=& \cfrac{CM_{k}}{N_{k}} \textrm{ for } 2\leq k \leq n ,\\
\hat{W}_{n-k} &=& 1+\cfrac{M_{k+1} -N_{k-1}}{N_{k}} \textrm{ for } 2\leq k \leq n-1.
\end{eqnarray*}
Then, we set
\[ U_{k} := \cfrac{M_{k}}{N_{k}} \textrm{ and } V_{k} := \cfrac{M_{k+1}}{N_{k}}  \; \textrm{ for } 2\leq k \leq n ,\]
and $V_{1}=\frac{1-2A_2}{A_2}$ getting the formulae
\begin{eqnarray*}
\hat{W}_{0} &= &V_{n-1},\\
W_{n-k}&=& CU_{k} \textrm{ for } 2\leq k \leq n ,\\
\hat{W}_{n-k} &=& 1+ V_{k} - \cfrac{U_{k}}{V_{k-1}} \textrm{ for } 2\leq k \leq n-1.
\end{eqnarray*}
Hence, the $2n-2$ independent parameters $C, A_{3}, \cdots , A_{n}, B_{2}, \cdots, B_{n}$ provide the $2n-2$ independent parameters $C,U_{2}, \ldots , U_{n}, V_{2}, \ldots , V_{n-1}$. Thus, it is sufficient to show that the rational fractions $W_{0}, \ldots , W_{n-2},\hat{W}_{0}, \ldots , \hat{W}_{n-2}\in \mathbb{Q}(A_2)(C,U_{2}, U_{3}, \ldots , U_{n},V_{2}, V_{3}, \ldots , V_{n-1})$ are algebraically independent over $\mathbb{Q}(A_2)$.\\

Suppose that there exists an irreducible polynomial $R \in \mathbb{Q}(A_2)[X_{1}, \ldots, X_{2n-2}]$ such that

\[ R \left(\hat{W}_{0}, \ldots , \hat{W}_{n-2},W_{0}, \ldots , W_{n-2} \right)=0.\]
 
Specializing $C$ in $1$, we obtain
\begin{equation}\label{C1} R\left(  V_{n-1}, V_{n-1} +1 - \cfrac{U_{n-1}}{V_{n-2}}, \ldots , V_{2} +1 - \cfrac{U_{2}}{V_1}  ,U_{n}, \ldots , U_{2} \right) =0\end{equation}
where the $2n-3$ last rational fractions generate the field $\mathbb{Q}(A_2)(U_{2}, \ldots , U_{n}, V_{2}, \ldots , V_{n-1})$ over $\mathbb{Q}(A_2)$. Therefore, they are algebraically independent. We investigate their relation with the first coordinate, that will provide information on $R$. Observe that for $2\leq k \leq n-1$,
\[\hat{W}_{n-k} = 1+ V_{k} - \cfrac{U_{k}}{V_{k-1}} \]
provide the relation
\[ V_{k}= \hat{W}_{n-k} -1 + \cfrac{{W}_{n-k} }{V_{k-1}}.\]
Since $\hat{W}_{0}= V_{n-1}$, we can compute by finite induction
\[ \hat{W}_{0}= V_{n-1} = (\hat{W}_1-1) + \frac{W_1}{V_{n-2}} = f_0 + \underset{k=1}{\overset{n-2}{\mathrm K}}\frac{e_k}{f_k}\]
where 
\[\left\{\begin{array}{clr} 
e_k &= W_k &\textrm{ for } 1\leq k \leq n-2\\
f_k &= \hat{W}_{k+1}-1 &\textrm{ for } 0 \leq k \leq n-3\\ 
f_{n-2} &= V_1 = \cfrac{1-2A_2}{A_2}&
\end{array}\right.\]
and \[f_0 + \underset{k=1}{\overset{n-2}{\mathrm K}}\frac{e_k}{f_k} = f_0 + \cfrac{e_1}{f_1 + \cfrac{e_2}{f_2 + \cfrac{\ddots}{f_{n-2}}}}\] is Gauss' notation for a (finite) generalized continued fraction. Denote by $\left(\cfrac{E_k}{F_k}\right)_{k=0}^{n-2}$ the finite sequence of its convergents. \\

We set 
\[ \tilde{R}=  F_{n-2}\hat{W}_0-E_{n-2}\]
where $F_{n-2}$ and $E_{n-2}$ are seen as polynomials in $\mathbb{Q}(A_2)[W_0, \ldots , W_{n-2}, \hat{W}_0, \ldots , \hat{W}_{n-2}].$ Note that $F_{n-2}$ and $E_{n-2}$ do not depend on $\hat{W}_0$ since none of the $(e_k)_{1\leq k \leq n-2}$ and $(f_k)_{0\leq k \leq n-2}$ do. Hence, $\tilde{R}$ is a polynomial of degree $1$ with respect to $\hat{W}_0$. Writing the Euclidean division of $R$ by $\tilde{R}$ in $\mathbb{Q}(A_2, \hat{W}_1, \ldots, \hat{W}_{n-2}, W_0, \ldots , W_{n-2})[\hat{W}_0]$ we get
\[ R = \tilde{R}Q + P\]

with $\deg_{\hat{W}_0}(P) = 0$. Hence $P$ can be seen as a polynomial in the $2n-3$ variables $\hat{W}_1, \ldots , \hat{W}_{n-2}, W_0, \ldots , W_{n-2}$ over $\mathbb{Q}(A_2)$. The latter are algebraically independent over $\mathbb{Q}(A_2)$ because their specializations at $C=1$ are. We deduce that $P=0$, and by irreducibility of $R$, the polynomial $Q$ is a constant:

\[ R = \alpha\left(F_{n-2} \hat{W}_0 - E_{n-2}\right) \]
with $\alpha \in \mathbb{Q}(A_2)$.\\

Specializing $C$ in $0$, we obtain
\[ R\left(  V_{n-1}, V_{n-1} +1 - \cfrac{U_{n-1}}{V_{n-2}}, \ldots , V_{2} +1 - \cfrac{U_{2}}{V_1}  ,0, \ldots , 0 \right) =0\]
where the $n-1$ non zero rational fractions generate the field $\mathbb{Q}(V_1)(U_3, \ldots, U_{n-1})(V_{n-1}, V_{n-2}, \ldots , V_{2}, U_2)$ over $\mathbb{Q}(V_1)(U_3, \ldots, U_{n-1})$. Therefore, they are algebraically independent over $\mathbb{Q}(A_2)=\mathbb{Q}(V_1)$. We deduce that the constant monomial of $R$ seen in $\mathbb{Q}(A_2, \hat{W}_{0}, \ldots , \hat{W}_{n-2})[{W}_{0}, \ldots , {W}_{n-2}]$ should be zero.\\

We now compute the constant monomial of $F_{n-2}\hat{W}_0-E_{n-2}$ seen in $\mathbb{Q}(A_2,\hat{W}_{0}, \ldots , \hat{W}_{n-2})[{W}_{0}, \ldots , {W}_{n-2}]$. We use the classical recurrence formulae for the convergents
\[E_{k+1} = e_{k+1} E_k + f_{k+1} E_{k-1} \textrm{ and } F_{k+1} = e_{k+1} F_k + f_{k+1} F_{k-1}\]
to compute the constant term of $E_{n-2}$ and $F_{n-2}$ to be $$\prod_{k=0}^{n-2} f_k \textrm{ and }\prod_{k=1}^{n-2} f_k$$ respectively. Thus the constant monomial of $F_{n-2}\hat{W}_0-E_{n-2}$ seen in $\mathbb{Q}(A_2,\hat{W}_{0}, \ldots , \hat{W}_{n-2})[{W}_{0}, \ldots , {W}_{n-2}]$ is
\[ \left(\prod_{k=1}^{n-2} f_k \right)\hat{W}_0 - \prod_{k=0}^{n-2} f_k = (\hat{W}_0-\hat{W}_1 +1)\cfrac{1-2A_2}{A_2}\prod_{k=1}^{n-3} (\hat{W}_{k+1}-1).\]
The fact that $\hat{W}_{k+1}\neq1$ and $\hat{W}_0+1 \neq \hat{W}_1 $ induces that this constant monomial is non zero. Hence $\alpha$ and  $R$ are zero.\\

This proves the algebraic independence of the $2n$ exponents. \qed

 \bibliographystyle{plain}
 \bibliography{Biblio}

 \end{document}